\documentclass[12pt,reqno]{amsart}
\renewcommand{\baselinestretch}{1.26}
\usepackage{epsfig}
\usepackage{rotating}
\usepackage{amssymb}
\usepackage{amsthm}
\usepackage{xcolor}
\usepackage{graphicx}
\usepackage{verbatim} 
\makeatletter
\newcommand{\singlespacing}{\let\CS=\@currsize\renewcommand{\baselinestretch}{1}\tiny\CS}
\begin{document}

\parskip = 10pt
\def \qed {\hfill \vrule height7pt width 5pt depth 0pt}
\newcommand{\ve}[1]{\mbox{\boldmath$#1$}}
\newcommand{\IR}{\mbox{$I\!\!R$}}
\newcommand{\1}{\Rightarrow}
\newcommand{\bs}{\baselineskip}
\newcommand{\esp}{\end{sloppypar}}
\newcommand{\be}{\begin{equation}}
\newcommand{\ee}{\end{equation}}
\newcommand{\beanno}{\begin{eqnarray*}}
\newcommand{\inp}[2]{\left( {#1} ,\,{#2} \right)}
\newcommand{\eeanno}{\end{eqnarray*}}
\newcommand{\bea}{\begin{eqnarray}}
\newcommand{\eea}{\end{eqnarray}}
\newcommand{\ba}{\begin{array}}
\newcommand{\ea}{\end{array}}
\newcommand{\nno}{\nonumber}
\newcommand{\dou}{\partial}
\newcommand{\bc}{\begin{center}}
\newcommand{\ec}{\end{center}}
\newcommand{\2}{\subseteq}
\newcommand{\cl}{\centerline}
\newcommand{\boxeq}[2]

\def\refhg{\hangindent=20pt\hangafter=1}
\def\refmark{\par\vskip 2.50mm\noindent\refhg}

\def \qed {\hfill \vrule height7pt width 5pt depth 0pt}

\def\lam{\lambda }
\def\Lam{\Lambda}
\def\lab{\label }
\def\iomui{[0,1/ \mu)}
\def\komu{\frac{k}{\mu}}
\def\kmu{\frac{k}{\mu}}
\def\Komu{\frac{K}{\mu}}
\def\alst{\alpha\sim t}
\def\tal{{\tau}^{\al}}
\def\abal{\mid{\alpha}\mid}
\def\la{\langle}
\def\ra{\rangle}
\def\rar{\rightarrow}
\def\Rar{\Rightarrow}
\def\ho{\hat{\omega}}
\def\to{\tilde{\omega}}
\def\o{\omega}

\def\cA{{\cal A}}
\def\cB{{\cal B}}
\def\cC{{\cal C}}
\def\cD{{\cal D}}
\def\cE{{\cal E}}
\def\cF{{\cal F}}
\def\cG{{\cal G}}
\def\cH{{\cal H}}
\def\cI{{\cal I}}
\def\cJ{{\cal J}}
\def\cK{{\cal K}}
\def\cL{{\cal L}}
\def\cM{{\cal M}}
\def\cN{{\cal N}}
\def\cO{{\cal O}}
\def\cP{{\cal P}}
\def\cQ{{\cal Q}}
\def\cR{{\mathcal R}}
\def\cS{{\cal S}}
\def\cT{{\cal T}}
\def\cU{{\cal U}}
\def\cV{{\cal V}}
\def\cW{{\cal W}}
\def\cX{{\cal X}}
\def\cY{{\cal Y}}
\def\cZ{{\cal Z}}

\def\1N{\frac{1}{n}}
\def\lb{\lbrace}
\def\rb{\rbrace}
\def\blp{\bigl (}
\def\brp{\bigr )}
\def\blb{\bigl \lbrace}
\def\brb{\bigr \rbrace}
\def\bls{\bigl [ }
\def\brs{\bigr ] }
\def\Blp{\Bigl (}
\def\Brp{\Bigr )}
\def\Blb{\Bigl \lbrace}
\def\Brb{\Bigr \rbrace}
\def\BLB{\Biggl \lbrace}
\def\BRB{\Biggr \rbrace}
\def\Bls{\Bigl [ }
\def\Brs{\Bigl ] }
\def\BLS{\Biggl [ }
\def\BRS{\Biggl ] }
\def\Skor{D([0, 1], \cS' (R^d))}
\def\mlnm{(m\lambda+n\mu)}
\def\ds{\displaystyle}
\def\what{\widehat}
\def\wtilde{\widetilde}

\def\limn{\lim_{n \rightarrow \infty}}
\def\limiN{\liminf_{n\rightarrow \infty}}
\def\stN{\sum_{t=1}^n}
\def\stn{\sum_{t=1}^n}
\def\sjp{\sum_{j=1}^p}
\def\skp{\sum_{k=1}^p}
\def\skii{\sum_{k=0}^\infty}
\def\slii{\sum_{l=0}^\infty}

\def\ptN{\prod_{t=1}^n}
\def\ptn{\prod_{t=1}^n}
\def\wh{\widehat}
\def\wt{\widetilde}
\def\a012{a_0 + a_1 t + a_2 t^2}
\def\0a012{a_0^0 + a_1^0 t + a_2^0 t^2}
\def\0t012{\theta_0^0 + \theta_1^0 t + \theta_2^0 t^2}
\def\t012{\theta_0 + \theta_1 t + \theta_2 t^2}
\def\blfootnote{\xdef\@thefnmark{}\@footnotetext} 
\theoremstyle{plain}
\newtheorem{theorem}{\bf \sc Theorem}[section]
\theoremstyle{definition}
\newtheorem{remark}{\bf Remark}
\theoremstyle{plain}
\newtheorem{lemma}{\bf Lemma}
\newtheorem{assum}{\bf Assumption}
\theoremstyle{definition}
\newtheorem{pth}{\bf Proof of Theorem}[]
\newtheorem{plm}{\bf Proof of Lemma}[]

\newcommand{\dotx}{\dot{X}}
\newcommand{\ddotx}{\ddot{X}}

\title[Fundamental Frequency Model]{\sc Estimating the fundamental frequency  using modified Newton-Raphson algorithm}
\author{Swagata Nandi$^1$}
\author{Debasis Kundu$^2$}
\address{$^1$Theoretical Statistics and Mathematics Unit, Indian Statistical Institute, 7, S.J.S. Sansanwal Marg, New Delhi - 110016,
India, nandi@isid.ac.in}
\address{$^2$Department of Mathematics and Statistics, Indian Institute of
Technology Kanpur, Pin 208016, India, kundu@iitk.ac.in} 

\keywords{Fundamental frequency model;  Approximate least squares estimator; Modified Newton-Raphson algorithm; Super efficient estimator;}
\subjclass[2000]{62J02; 62E20; 62C05}

\begin{abstract}
In this paper, we propose a modified Newton-Raphson algorithm to estimate the frequency parameter in the fundamental frequency model 
in presence of an additive stationary error.  The proposed estimator is super efficient in nature in the sense that its asymptotic 
variance is less than the asymptotic variance of the least squares estimator.  With a proper step factor modification, the 
proposed modified Newton-Raphson algorithm produces an estimator with the rate $O_p(n^{-\frac{3}{2}})$, the same rate as the least squares estimator.  Numerical experiments are performed for different sample sizes, different error variances and for different models.  For 
illustrative purposes, two real data sets are analyzed using the fundamental frequency model and the estimators are obtained using the proposed algorithm.  It is observed the model and the proposed algorithm work quite well in both cases. 
\end{abstract}

\maketitle

\section{Introduction}   \label{sec1}

In this paper, we consider the problem of estimating the frequency present in the following fundamental frequency model:
\begin{equation}
 y(t) = \sum_{j=1}^p [A_j \cos(j \lambda t) + B_j \sin(j \lambda t) ] + e(t), \hspace{.15in}   \label{fundamental_modeleq1}
t=1,\ldots, n
\end{equation}
Here $y(t)$ is the observed signal at time point $t$; $A_k, B_k \in \mathbb{R}$ are unknown amplitudes and none of them are 
identically equal to zero; $\ds 0 < \lambda < \pi/p$,  is the fundamental 
frequency;  the number of components $p$ is assumed to be known.  The sequence of error random variables $\{e(t)\}$ is from a stationary linear process and  satisfies the following assumption.
\begin{assum}  \label{ass1}
The sequence $\{e(t)\}$ has the following representation:
\begin{equation}
e(t) = \sum_{k=0}^\infty a(k) \epsilon(t -k), \hspace{.25in} \sum_{k=0}^\infty | a(k) | < 
\infty,   \label{assum_eq}
\end{equation}
where $\{\epsilon(t)\}$ is a sequence of independent and identically distributed (i.i.d.) random variables with mean zero and finite variance 
$\sigma^2 > 0$.  The arbitrary real-valued sequence $\{a(k)\}$ is absolutely summable.
\end{assum} 
Assumption \ref{ass1} is a standard assumption of a weakly stationary linear process.  Any stationary ARMA process satisfies Assumption \ref{ass1} and can be expressed as \eqref{assum_eq}. 
The fundamental frequency model \eqref{fundamental_modeleq1} is a very useful model for periodic signals when harmonics of a fundamental frequency are present.  The model has applications in a variety of fields and is a particular case of the usual sinusoidal model 
\begin{equation}
 y(t) = \sum_{j=1}^p [A_j \cos(\lambda_j t) + B_j \sin(\lambda_j t) ] + e(t), \hspace{.15in}   \label{sinusoidal_model}
t=1,\ldots, n.
\end{equation}
The model \eqref{fundamental_modeleq1} is a particular case of the model \eqref{sinusoidal_model} with a restriction in model parameters; the frequency of the $j$th component of the sinusoidal model $\lambda_j = j \lambda$.  When frequencies are at $\lambda, 2\lambda, \ldots, p\lambda$ instead of arbitrary $\lambda_j \in (0, \pi)$, $j=1,\ldots,p$, these are termed as harmonics of $\lambda$.   The presence of an exact periodicity is a convenient approximation, but many real life phenomena can be described using model \eqref{fundamental_modeleq1}.   There are many non-stationary signals like speech, human circadian system where the data indicate the presence of harmonics of a fundamental frequency.   In such cases, it is more convenient to use model \eqref{fundamental_modeleq1} than \eqref{sinusoidal_model} because model \eqref{fundamental_modeleq1} has one non-linear parameter as compared to $p$ in model \eqref{sinusoidal_model}.

In the literature, many authors considered the following model instead of model \eqref{fundamental_modeleq1},
\begin{equation}
y(t) = \sum_{j=1}^p \rho_j \cos(tj\lambda - \phi_j) + e(t), \label{differ_para_funda_model}
\end{equation}
where $\rho_j$'s are amplitudes, $\lambda$ is the fundamental frequency and $\phi_j$'s are phases and $\rho_j >0$, 
$\lambda \in (0, \pi/p)$ and $\phi_j \in (0,\pi)$, $j=1,\ldots, p$. The sequence $\{e(t)\}$ is the error component.  Note that model \eqref{differ_para_funda_model} is same as model \eqref{fundamental_modeleq1} with a different parameterization.  In this case, $A_j=\rho_j \cos(\phi_j)$ and $B_j = -\rho_j \sin(\phi_j)$.

We are mainly interested to estimate the fundamental frequency present in model \eqref{fundamental_modeleq1} under assumption \ref{ass1}.   The problem was originally proposed by Quinn and Thomson \cite{QT:1991} and they proposed a weighted least squares approach to estimate the unknown parameters.  This is  basically an approximate generalized least squares criterion.  Since then an extensive 
amount of work has been done dealing with different aspects of this model.  
Abatzoglou et al. \cite{AMH:1991} considered the total least squares estimators approach in estimating the parameters of the above model.
Nandi and Kundu \cite{NK:2004, NK:2003} studied the asymptotic properties of the least squares estimator of the unknown parameters of the 
model \eqref{differ_para_funda_model} under Assumption \ref{ass1}.  Cristensen et al. \cite{CJJ:2007} proposed joint estimation of fundamental frequency and number of harmonics based on MUSIC criterion.  
Recently, Nielsen et al. \cite{NJJCJ:2017}
provided a computationally efficient estimator of the unknown parameters of the model \eqref{fundamental_modeleq1}.
A more general model with presence of multiple fundamental frequencies has been considered by Christensen et al. \cite{CHJJ:2011} and Zhou \cite{Zhou:2013}.  A further generalized model where fundamental frequencies appear in clusters has been proposed by Nandi and Kundu \cite{NK:2006a}.

It is a well known fact that even for the usual sinusoidal model \eqref{sinusoidal_model}, the Newton-Raphson (NR) algorithm does not work well, see for example Bresler and Macovski \cite{BM:1986}.  In many situation the NR algorithm does not converge or converges to a local minimum.  In this paper, we have modified the Newton-Raphson algorithm by reducing the step factor in the NR algorithm applied to an equivalent criterion function of the approximate least squares estimator.   We have proved that the estimator obtained from the modified NR 
algorithm has the same rate of convergence as the LSEs.  Moreover, the asymptotic variance of the modified NR (MNR) estimate is 
one fourth of the asymptotic variance of the least squares estimator.  Our approach is similar to the approach adopted by Kundu et al. \cite{KBNB:2011} or Bian et al. 
\cite{BPXLL:2013}.  Kundu et al. \cite{KBNB:2011} considered the model (\ref{sinusoidal_model}) in presence of additive noise, and 
Bian et al. \cite{BPXLL:2013} also considered the same model in presence of multiplicative and additive error.  In both the cases the 
authors obtained super efficient estimators of $p$ frequencies in a sequential approach.  But in this case the model has one fundamental 
frequency and $p$ harmonics, hence the sequential procedure is not possible.  The fundamental frequency and the harmonics need to be 
estimated simultaneously.  This is the main difference of the present manuscript with the existing work.

Model \eqref{fundamental_modeleq1} is an important model in analyzing periodic data and can be useful in situations where periodic signals are observed with an inherent fundamental frequency.  Baldwin and Thomson \cite{BT:1978} used model \eqref{fundamental_modeleq1} to describe the visual observations of S.Carinae, a variable star in the sky of the Southern Hemisphere.  Greenhouse et al. \cite{GKT:1987} proposed the use of higher-order harmonic terms of one or more fundamentals and ARMA
processes for the errors for fitting biological rhythms (human core body temperature data).  For illustrative purposes we use model 
\eqref{fundamental_modeleq1} to analyze two  speech data sets and estimate the parameters using the proposed algorithm.  It is observed that 
the model and the proposed algorithm work quite satisfactorily in both the cases.  

The rest of the paper is organized as follows.  In section \ref{sec2}, the least squares and the approximate least squares criteria for the fundamental frequency model are described.  In section \ref{sec3}, we propose the MNR algorithm and state the main result of the paper. In section \ref{sec4}, we carry out  numerical experiments based on simulation.  Two real speech data sets are analyzed for illustrative purposes in section \ref{sec5}, and in final section, we summarize the results and directions for future work.
  
\section{Estimation of Unknown Parameters} \label{sec2}

In matrix notation, model \eqref{fundamental_modeleq1} can be written as 
\begin{equation}
{\bf Y} = {\bf X}(\lambda) {\ve \theta} + {\bf e},
\end{equation}
where ${\bf Y} = (y(1), \ldots, y(n))^T$,  ${\bf e} = (e(1), \ldots, e(n))^T$,  ${\ve \theta} = (A_1, B_1, \ldots, A_p, B_p)^T$, 
${\bf X}(\lambda) = ({\bf X}_1, \ldots, {\bf X}_p)$ and 
$$
{\bf X}_j = \left[\begin{matrix}
\cos(j \lambda) & \sin(j \lambda) \\
\cos(2 j \lambda) & \sin(2 j \lambda) \\
\vdots & \vdots \\
\cos(n j \lambda) & \sin(n j \lambda) \\
\end{matrix}\right].
$$

The matrix ${\bf X}_j = {\bf X}_j(\lambda)$, but we do not make it explicit.
The least squares criterion minimizes 
\begin{equation}
Q({\ve \theta}, \lambda) = ({\bf Y} - {\bf X}(\lambda){\ve \theta})^T ({\bf Y} - {\bf X}(\lambda){\ve\theta}).
\end{equation}
For a given $\lambda$, $Q({\ve\theta}, \lambda)$ is minimized at $\what{\ve\theta}(\lambda) = ({\bf X}(\lambda)^T  {\bf X}
(\lambda))^{-1} {\bf X}(\lambda)^T {\bf Y}$.  Then,
\begin{eqnarray*}
Q(\what{\ve\theta}, \lambda) &=& \left({\bf Y} - {\bf X}(\lambda)\what{\ve\theta}(\lambda)\right)^T  \left({\bf Y} - {\bf X}(\lambda
)\what{\ve\theta}(\lambda) \right) \\
&=& {\bf Y}^T  {\bf Y} -  {\bf Y}^T {\bf X}(\lambda)  \left({\bf X}(\lambda)^T ({\bf X}(\lambda) \right)^{-1} {\bf X}(\lambda)^T {\bf Y}.
\end{eqnarray*}
Therefore, minimizing $Q(\what{\ve\theta}, \lambda)$ with respect to $\lambda$ is equivalent to maximizing 
$$
{\bf Y}^T {\bf X}(\lambda) \left({\bf X}(\lambda)^T ({\bf X}(\lambda) \right)^{-1} {\bf X}(\lambda)^T  {\bf Y},
$$
with respect to $\lambda$.  This quantity is asymptotically equivalent to (see Nandi and Kundu \cite{NK:2003})
\begin{equation}
Q_N(\lambda) = \sum_{j=1}^p \left| \frac{1}{n} \sum_{t=1}^n y(t) e^{itj\lambda} \right|^2.   \label{quinn_criterion}
\end{equation}
On the other hand, ${\bf Y}^T {\bf X}_j ({\bf X}_j^T {\bf X}_j)^{-1} {\bf X}_j^T {\bf Y}$ and 
$\left| \left (\sum_{t=1}^n y(t) e^{itj\lambda} \right )/n \right|^2$
 are asymptotically equivalent.  Hence, the criterion is based on the maximization of 
\begin{equation}
g(\lambda) = \sum_{j=1}^p \left[ {\bf Y}^T {\bf X}_j ({\bf X}_j^T {\bf X}_j)^{-1} {\bf X}_j^T {\bf Y} \right],
\end{equation}
with respect to $\lambda$.   Write ${\bf Y}^T {\bf X}_j ({\bf X}_j^T {\bf X}_j)^{-1} 
{\bf X}_j^T {\bf Y} = R_j(\lambda)$, then 
\begin{equation}
\what{\lambda} = \arg\max_{\lambda} g(\lambda) = \arg\max_{\lambda} \sum_{j=1}^p R_j(\lambda). \label{lambdahat_eq}
\end{equation}
Note that for large $n$, $\left ({\bf X}_j^T {\bf X}_k \right )/n = {\bf 0}$, for $j\ne k$.  Hence, 
$$
Q(\what{\ve\theta}, \lambda)  = {\bf Y}^T{\bf Y} - \frac{1}{n} \sum_{j=1}^p {\bf Y}^T {\bf X}_j^T ({\bf X}_j^T {\bf X}_j)^{-1} {\bf X}_j^T {\bf Y}.
$$ 
Once $\what{\lambda}$ is estimated using \eqref{lambdahat_eq}, the linear parameters are either estimated as least squares estimators,
$$
\left(\begin{matrix}
\what{A}_j \\ \what{B}_j \\
\end{matrix} \right) = ({\bf X}_j(\what{\lambda})^T  {\bf X}_j(\what{\lambda}))^{-1} {\bf X}_j(\what{\lambda})^T {\bf Y}.
$$
or as approximated least squares estimators, given as follows:
\be
\widetilde{A}_j = \frac{2}{n} \sum_{t=1}^n y(t) \cos(j\lambda t), \;\;\; \widetilde{B}_j = \frac{2}{n} \sum_{t=1}^n y(t) \sin(j\lambda t).
\ee
The estimator of $\lambda$ defined in \eqref{lambdahat_eq} is nothing but the approximate least squares estimator (ALSE) of $\lambda$ which has been studied extensively in the literature.

The asymptotic distribution of the least squares estimators and approximate least squares estimators of the unknown parameters of 
model \eqref{fundamental_modeleq1} under assumption \ref{ass1} are obtained by Nandi and Kundu \cite{NK:2003}.  In fact, Nandi and Kundu \cite{NK:2003} studied model \eqref{differ_para_funda_model} and observed that the asymptotic distribution of LSEs and ALSEs are same.  Under assumption \ref{ass1}, the asymptotic distribution of $\what{\lambda}$, the LSE of $\lambda $ is as follows:
\begin{equation}
n^{3/2} (\what{\lambda} - \lambda)  \stackrel{d}{\longrightarrow} \mathcal{N} \left(0, \frac{24\sigma^2 \delta_G}{{\beta^*}^2} \right),
  \label{asym_lse}
\end{equation} 
where $\ds \beta^* = \sum_{j=1}^p j^2 (A_j^2 + B_j^2)$, $\ds \delta_G = \sum_{j=1}^p j^2 (A_j^2 + B_j^2) c(j)$ and 
 $\ds c(j)= \left| \sum_{k=0}^\infty a(k)e^{-ijk\lambda} \right|^2$.  The notation $\stackrel{d}{\longrightarrow}$ means convergence in distribution and $\mathcal{N}(a,b)$ denotes the Gaussian distribution with mean $a$ and variance $b$.  

\section{Modified Newton-Raphson Algorithm}  \label{sec3}

We first describe the standard NR algorithm  in case of $\ds g(\lambda) = \sum_{j=1}^p R_j(\lambda)$, before proceeding further.  Let $\what{\lambda}_1$ be the initial estimate of $\lambda$ and $\what{\lambda}_k$ be the  estimate at the $k$th iteration.  Then, the NR estimate at the $(k+1)$th iteration is obtained as 
\begin{equation}
\what{\lambda}_{k+1} = \what{\lambda}_{k} - \frac{g'(\what{\lambda}_k)}{g''(\what{\lambda}_k)},  \label{nr_eq}
\end{equation}
where $g'(\what{\lambda}_k)$ and $g''(\what{\lambda}_k)$ are first and second order derivatives of $g(.)$ evaluated at $\what{\lambda}_k$, respectively.

The standard NR algorithm \eqref{nr_eq} is modified by reducing the step factor as follows:
\begin{equation}
\what{\lambda}_{k+1} = \what{\lambda}_{k} - \frac{1}{4}\frac{g'(\what{\lambda}_k)}{g''(\what{\lambda}_k)}.   \label{modified_nr_eq}
\end{equation}
A smaller step factor prevents the algorithm to diverge.  At a particular iteration, if the estimator is close enough to the global minimum, then a comparatively large correction factor may shift the estimate  far away from the global minimum.  Therefore, with the 
proper choice of the initial guess and with the correct step factor, the algorithm provides a super efficient estimator of 
$\lambda$.  We need the following theorem for further development.

\begin{theorem}  \label{theorem1}
Let $\wt{\lambda}_0$ be a consistent estimator of $\lambda$ and $\wt{\lambda}_0 - \lambda = O_p(n^{-1-\delta})$, $\delta \in (0, \frac{1}{2}]$.  Suppose $\wt{\lambda}_0$ is updated as $\ds \wt{\lambda} = \wt{\lambda}_0 - \frac{1}{4}\frac{g'(\wt{\lambda}_0)}{g''(\wt{\lambda}_0)}$, then
\begin{enumerate}
\item[{(a)}] $\ds \wt{\lambda} - \lambda  = O_p(n^{-1-3\delta})$ if $\delta \le \frac{1}{6}$. \\

\item[{(b)}] $\ds n^{3/2}(\wt{\lambda} - \lambda) \rightarrow \mathcal{N} \left(0, \frac{6 \sigma^2 \delta_G}{{\beta^*}^2} \right)$, if $\delta > \frac{1}{6}$, 
\end{enumerate}
where $\beta^*$ and $\delta_G$ are same as defined in the previous section.  
\end{theorem}

\noindent {\sc Proof:} See in the Appendix.   \qed

It should be mentioned that the step factor $1/4$ is not arbitrary.  The motivation to take the step factor as 
$1/4$, and the reason that the algorithm provides a super efficient estimator come from Theorem \ref{theorem1}.  It can be seen from the last two 
equalities in (\ref{proof_eq7}), and that is the key step why the method works.  Theorem \ref{theorem1} immediately shows how the improvement in order can be made from $\wt{\lambda}$ to $\wh{\lambda}$.  But now
if we take the step factor in (\ref{modified_nr_eq}) anything smaller than $1/4$, say, $1/8$, then it can be shown along the same way 
that no improvement in order can be made.  It means that the step is too small.  
Interestingly, if we take the step factor in (\ref{modified_nr_eq}) anything greater than $1/4$, say, $1/2$, then also it can be shown that
the improvement in order cannot be made.  In this case the step is too large and it crosses 
the target and moves to the other direction.  In fact, that is the reason the standard Newton-Raphson algorithm does not work
in this case.

This theorem states that if we start from a reasonably good initial estimator, then the MNR algorithm produces estimator with the same convergence rate as that of the LSE of $\lambda$.  Moreover, the asymptotic variance of the proposed estimator of the fundamental frequency 
is one fourth of the  asymptotic variance of the LSE.  The argument maximum of the periodogram function over Fourier frequencies $\frac{2\pi k}{n}, k=1,\ldots, \left[\frac{n}{2}\right]$ , as an estimator of the frequency has a convergence rate $O_p(n^{-1})$.  We use this estimator as the starting value of the algorithm implemented with a subset of the observed data vector of size $n$ using similar idea of Kundu et al. \cite{KBNB:2011} and Bian et al. \cite{BPXLL:2013}.  The subset is selected in such a way that the dependence structure present in the data is not destroyed, that is, a subset of predefined size of consecutive points is selected as a starting point.  
The details have been illustrated in the data analyses section.    

\noindent {\bf Algorithm:}

\begin{enumerate}
\item  Obtain the argument maximum of the periodogram function $I(\lambda)$ over Fourier frequencies and denote it as $\widetilde{\lambda}_{0}$.  Then 
$\widetilde{\lambda}_{0}= O_p(n^{-1})$.

\item At $k=1$, take $n_1 = n^{6/7}$ and calculate $\wt{\lambda}_1$ as 
\be
\wtilde{\lambda}_{1} = \wtilde{\lambda}_{0} - \frac{1}{4}\frac{g'_{n_1}(\wtilde{\lambda}_0)}{g''_{n_1}(\wtilde{\lambda}_0)}.
\ee
where $g'_{n_1}(\wtilde{\lambda}_0)$ and $g''_{n_1}(\wtilde{\lambda}_0)$ are same as $g'(.)$ and $g''(.)$ evaluated at $\wtilde{\lambda}_0$ with a sub-sample of size $n_1$.  Note that $\widetilde{\lambda}_{0} - \lambda = O_p(n^{-1})$ and $n_1 = n^{6/7}$, so $n=n_1^{-7/6}$ .  Therefore, 
$\widetilde{\lambda}_{0} - \lambda = O_p(n^{-1}) = O_p(n_1^{-1-\frac{1}{6}})$ and applying part (a) of theorem \ref{theorem1},  we have 
$
\widetilde{\lambda}_{1} - \lambda = O_p(n_1^{-1-\frac{1}{2}}) = O_p(n^{-\frac{3}{2} \times \frac{6}{7}}) =  O_p(n^{-\frac{9}{7}}) = O_p(n^{-1-\frac{2}{7}})
$ with $\delta = \frac{2}{7}$. 

\item As $\widetilde{\lambda}_{1} - \lambda = O_p(n^{-1-\frac{2}{7}})$, $\delta =\frac{2}{7} > \frac{1}{6}$,  we can apply part (b) of theorem \ref{theorem1}.   Take $n_{k+1}=n$, and repeat
 \begin{equation}
\wt{\lambda}_{k+1} = \wt{\lambda}_{k} - \frac{1}{4}\frac{g'_{n_k+1}(\wt{\lambda}_k)}{g''_{n_k+1}(\wt{\lambda}_k)}, \;\;\;k=1,2, \ldots   \label{final_alg_eq}
\end{equation}
until a suitable stopping criterion is satisfied.   \qed
\end{enumerate} 
Using theorem \ref{theorem1}, the algorithm implies that if at any steps, the estimator of $\lambda$ is of order $O_p(n^{-1-\delta})$, the updated estimator is of order $O_p(n^{-1-3\delta})$ if $\delta \le \frac{1}{6}$ and if $\delta > \frac{1}{6}$, the updated estimator is of same order as the LSE.  In addition, the asymptotic variance is four times less than the LSE, hence we call it a super efficient estimator.  In the proposed algorithm, we have started with a sub-sample of size $n^{6/7}$ of the original sample of size $n$.   The factor $\frac{6}{7}$  is not that important and not unique.  There are several other choices of $n_1$ to initiate the algorithm, for example, $n_1=n^{\frac{8}{9}}$ and $n_k=n$ for $k\ge 2$.

To obtain an estimator of order $O_p(n^{-1})$ is easy, but an estimator of order  $O_p(n^{-1-\delta})$, $\delta \in (0,\frac{1}{2}]$ is required to apply theorem \ref{theorem1}.  We have started the algorithm with a smaller number of observations to overcome this problem.  Varying sample size enables us to get estimator of order $O_p(n^{-1-\delta})$, for some $\delta \in (0,\frac{1}{2}]$.  With the particular choice of $n_1$, we can use all the available data points from second step onwards.   The proposed algorithm provides a super efficient estimator of the fundamental frequency from the relatively poor periodogram maximizer over the Fourier frequencies.  It is worth mentioning at this 
point that the initial estimator is not the ALSE and is not asymptotically equivalent to the LSE.  ALSE of $\lambda$ in case of fundamental frequency model maximizes the sum of $p$ periodogram functions at the harmonics without the constraints of Fourier frequencies (see Nandi and Kundu \cite{NK:2003}).

\section{Numerical Experiments}  \label{sec4}

In this section, numerical experimental results are presented based on Monte Carlo simulations to observe the performance of the proposed estimator.  We consider model \eqref{fundamental_modeleq1} with $p=4$ with two different sets of parameters as follows:
\beanno
\hbox{Model 1}: && A_1 = 5.0, \;\;\; A_2= 4.0, \;\;\; A_3 = 3.0, \;\;\; A_4 = 2.0, \\
&& B_1 = 3.0, \;\;\; B_2 = 2.5, \;\;\; B_3 = 2.25, \;\;\; B_4 = 2.0, \;\;\; \lambda = .25 \\
\hbox{Model 2}: && A_1 = 4.0, \;\;\; A_2= 3.0, \;\;\; A_3 = 2.0, \;\;\; A_4 = 1.0, \\
&& B_1 = 2.0, \;\;\; B_2 = 1.5, \;\;\; B_3 = 1.25, \;\;\; B_4 = 1.0, \;\;\;  \lambda = .3141.
\eeanno
The sequence of error random variables $\{e(t)\}$ is a moving average process of order one, $e(t) = .5 \epsilon(t-1) + \epsilon(t)$; $\epsilon(t)$ is a sequence of i.i.d. Gaussian random variables with mean zero and variance $\sigma^2$.  We consider different sample sizes, $n$ = 100, 200, 400, 500, 1000, and different error variance of $\{\epsilon(t)\}$, $\sigma^2$ =.01, .25, .75, and 1.0.  Note that, for the generated MA process $\{e(t)\}$, the variance is $1+\sigma^2$.   For the numerical experiments considered in this section, we assume that $p$ is known.  
 
In each case we generate a sample of size $n$ using the given model parameters and the error sequences.  We start the iteration with 
the given initial value.  The iterative process is terminated when the absolute difference between two consecutive iterates is less than $10^{-7}$ or the value of the objective function does not decrease.  In each case we report the average estimates and the variance of 
the estimates based on $5000$ replications.  The asymptotic variance of the proposed estimator as stated in theorem \ref{theorem1}(b) as well as the asymptotic variance of the LSE provided in \eqref{asym_lse} are also reported for comparison purposes.  The results for Model 1 are reported in Tables \ref{table1} and \ref{table2} and for Model 2, in Tables \ref{table3} and \ref{table4}.  It should be mentioned that when the
errors are i.i.d. and normally distributed then the asymptotic variance of the LSE is same as the Cramer-Rao lower bound.  But without 
any distributional assumptions on the error random variables, the Cramer-Rao lower bound cannot be computed.  

\begin{table}[t]
\caption{The average estimates, mean squared errors, asymptotic variances of LSEs and MNR estimates of the fundamental frequency in case Model 1 with correlated error when sample size $n=100$, $200$, $400$, $500$ and $1000$.}  
\label{table1}
\begin{tabular}{|c|c|c|c|c|}
\hline
\multicolumn{5}{c}{Sample Size N=100} \\
\hline  
$\sigma^2$  & Average & Variance & Asym. Var. (LSE) & Asym. Var. (MNR) \\
\hline
 .01   &    .252     &    8.07e-10   &     1.25e-9   &     3.13e-10     \\ \hline  
 .25       &   .252     &    1.73e-8   &     3.13e-8   &     7.84e-9 \\ \hline 
 .75       &    .252     &    5.13e-8   &     9.40e-8   &     2.35e-8   \\ \hline 
 1.0       &    .252        &    6.84e-8   &     1.25e-7   &     3.13e-8 \\ \hline 
\multicolumn{5}{c}{Sample Size N=200} \\
\hline
$\sigma^2$  & Average & Variance & Asym. Var. (LSE) & Asym. Var. (MNR) \\
\hline
.01   &   .250      &    9.93e-11   &     1.57e-10   &     3.92e-11   \\ \hline
.25      &    .250       &    2.37e-9   &     3.92e-9   &     9.80e-10   \\ \hline 
.75       &     .250        &    6.66e-9   &     1.18e-8   &     2.94e-9  \\ \hline 
1.0       &    .250       &    8.80e-9   &     1.57e-8   &     3.92e-9 \\ \hline 
\multicolumn{5}{c}{Sample Size N=400} \\
\hline
$\sigma^2$  & Average & Variance & Asym. Var. (LSE) & Asym. Var. (MNR) \\
\hline
.01   &  .250       &    1.85e-11   &     1.96e-11   &     4.90e-12   \\ \hline
.25       &    .250     &    3.93e-10   &     4.90e-10   &     1.22e-10   \\ \hline 
.75       &   .250      &    1.06e-9   &     1.47e-9   &     3.67e-10   \\ \hline 
1.0       &    .250       &    1.40e-9   &     1.96e-9   &     4.90-10  \\ \hline
\multicolumn{5}{c}{Sample Size N=500} \\
\hline
$\sigma^2$  & Average & Variance & Asym. Var. (LSE) & Asym. Var. (MNR) \\
\hline
.01   &   .250     &    9.51e-12   &     1.00e-11   &     2.51e-12    \\ \hline 
.25       &    .250      &    1.80e-10   &     2.51e-10   &     6.27e-11  \\ \hline  
.75       &    .250       &    5.17e-10   &     7.52e-10   &     1.88e-10  \\ \hline
 1.0       &    .250      &    6.85e-10   &     1.00e-9   &     2.51e-10   \\ \hline 
\multicolumn{5}{c}{Sample Size N=1000} \\
\hline
$\sigma^2$  & Average & Variance & Asym. Var. (LSE) & Asym. Var. (MNR) \\
\hline
.01   &   .250     &    6.06e-13   &     1.25e-12   &     3.13e-13   \\ \hline 
.25       &      .250       &    1.51e-11   &     3.13e-11   &     7.84e-12 \\ \hline 
.75       &     .250       &    4.93e-11   &     9.40e-11   &     2.35e-11   \\ \hline 
1.0       &   .250       &    6.08e-11   &     1.25e-10   &     3.13e-11    \\ \hline 
\end{tabular}
\end{table}

\begin{table}[t]
\caption{The average estimates, mean squared errors, asymptotic variances of LSEs and MNR estimates of the fundamental frequency in case Model 1 with i.i.d. error  when sample size $n=100$, $200$, $400$, $500$ and $1000$.}  
\label{table2}
\begin{tabular}{|c|c|c|c|c|}
\hline
\multicolumn{5}{c}{Sample Size N=100} \\
\hline  
$\sigma^2$  & Average & Variance & Asym. Var. (LSE) & Asym. Var. (MNR) \\
\hline
.01   &   .252       &    4.58e-10   &     6.36e-10   &     1.59e-10  \\ \hline
,25       &    .252      &    9.95e-9   &     1.59e-8   &     3.97e-9  \\ \hline
.75       &     .252       &    2.85e-8   &     4.77e-8   &     1.19e-8 \\ \hline  
1.0       &   .252        &    3.78e-8   &     6.36e-8   &     1.59e-8 \\ \hline
\multicolumn{5}{c}{Sample Size N=200} \\
\hline
$\sigma^2$  & Average & Variance & Asym. Var. (LSE) & Asym. Var. (MNR) \\
\hline
.01   &   .250      &    5.15e-11   &     7.95E-11   &     1.99e-11  \\ \hline
,25       &    .250       &    1.27e-9   &     1.99e-9   &     4.97e-10  \\ \hline
.75       &   .250       &    3.58e-9   &     5.96e-9   &     1.49e-9   \\ \hline  
1.0       &   .250       &    4.69e-9   &     7.95e-9   &     1.99e-9  \\ \hline
\multicolumn{5}{c}{Sample Size N=400} \\
\hline
$\sigma^2$  & Average & Variance & Asym. Var. (LSE) & Asym. Var. (MNR) \\
\hline
.01   &    .250     &    9.53e-12   &     9.93e-12   &     2.48e-12 \\ \hline
,25       &    .250       &    2.19e-10   &     2.48e-10   &     6.21e-11 \\ \hline
.75       &     .250       &    5.80e-10   &     7.45e-10   &     1.86e-10    \\ \hline  
1.0       &  .250      &    7.53e-10   &     9.93e-10   &     2.48e-10 \\ \hline
\multicolumn{5}{c}{Sample Size N=500} \\
\hline
$\sigma^2$  & Average & Variance & Asym. Var. (LSE) & Asym. Var. (MNR) \\
\hline
.01   &   .250      &    4.96e-12   &     5.09e-12   &     1.27e-12     \\ \hline
,25       &   .250       &    9.49e-11   &     1.27e-10   &     3.18e-11  \\ \hline
.75       &     .250       &    2.68e-10   &     3.81e-10   &     9.53e-11  \\ \hline  
1.0       & .250       &    3.54e-10   &     5.09e-10   &     1.27e-10  \\ \hline
\hline
\multicolumn{5}{c}{Sample Size N=1000} \\
\hline  
$\sigma^2$  & Average & Variance & Asym. Var. (LSE) & Asym. Var. (MNR) \\
\hline
.01   &   .250       &    3.12e-13   &     6.36e-13   &     1.59e-13  \\ \hline
,25       &    .250      &    6.28e-12   &     1.59e-11   &     3.97e-12  \\ \hline
.75       &     .250       &    2.06e-11   &     4.77e-11   &     1.19e-11 \\ \hline  
1.0       &   .250        &    2.68e-11   &     6.36e-11   &     1.59e-11 \\ \hline
\end{tabular}
\end{table}

\begin{table}[h]
\caption{The average estimates, mean squared errors, asymptotic variances of LSEs and MNR estimates of the fundamental frequency in case Model 2 correlated error when sample size $n=100$, $200$, $400$, $500$ and $1000$.}   \label{table3}
\begin{tabular}{|c|c|c|c|c|}
\hline
\multicolumn{5}{c}{Sample Size N=100} \\
\hline
$\sigma^2$  & Average & Variance & Asym. Var. (LSE) & Asym. Var. (MNR) \\
\hline
.01   &   .3148       &    1.67e-9   &     3.09e-9   &     7.73e-10  \\ \hline
.25       &    .3148       &    4.27e-8   &     7.73e-8   &     1.93e-8  \\ \hline 
.75       &     .3148       &    1.34e-7   &     2.32e-7   &     5.80e-8   \\ \hline  
1.0       &  .3148      &    1.80e-7   &     3.09e-7   &     7.73e-8    \\ \hline
\multicolumn{5}{c}{Sample Size N=200} \\
\hline
$\sigma^2$  & Average & Variance & Asym. Var. (LSE) & Asym. Var. (MNR) \\
\hline
.01   &   .3143       &    2.43e-10   &     3.86e-10   &     9.66e-11  \\ \hline
.25       &   .3143       &    6.11e-9   &     9.66e-9   &     2.42e-9   \\ \hline 
.75       &      .3143       &    1.74e-8   &     2.90e-8   &     7.25e-9  \\ \hline  
1.0       &    .3143      &    2.24e-8   &     3.86e-8   &     9.66e-9   \\ \hline
\multicolumn{5}{c}{Sample Size N=400} \\
\hline
$\sigma^2$  & Average & Variance & Asym. Var. (LSE) & Asym. Var. (MNR) \\
\hline
.01   &     .3141      &    3.59e-11   &     4.83e-11   &     1.21e-11  \\ \hline
.25       &    .3141       &    8.57e-10   &     1.21e-9   &     3.02e-10   \\ \hline 
.75       &      .3141       &    2.30e-9   &     3.62e-9   &     9.06e-10  \\ \hline  
1.0       &     .3141       &    3.03e-009   &     4.83e-9   &     1.21e-9  \\ \hline
\multicolumn{5}{c}{Sample Size N=500} \\
\hline
$\sigma^2$  & Average & Variance & Asym. Var. (LSE) & Asym. Var. (MNR) \\
\hline
.01   &    .3141      &    2.01e-11   &     2.47e-11   &     6.18e-12  \\ \hline
.25       &    .3141     &    4.58e-10   &     6.18e-10   &     1.55e-10   \\ \hline 
.75       &   .3141     &    1.22e-9   &     1.86e-9   &     4.64e-10  \\ \hline  
1.0       &   .3141       &    1.61e-9   &     2.47e-9   &     6.18e-10     \\ \hline
\multicolumn{5}{c}{Sample Size N=1000} \\
\hline
$\sigma^2$  & Average & Variance & Asym. Var. (LSE) & Asym. Var. (MNR) \\
\hline
.01   &   .3141       &    1.31e-12   &     3.09e-12   &     7.73e-13  \\ \hline
.25       &    .3141       &    3.88e-11   &     7.73e-11   &     1.93e-11  \\ \hline 
.75       &     .3141       &    1.01e-10   &     2.32e-10   &     5.80e-11   \\ \hline  
1.0       &  .3141      &    1.41e-10   &     3.09e-10   &     7.73e-11    \\ \hline
\end{tabular}
\end{table}

\begin{table}[h]
\caption{The average estimates, mean squared errors, asymptotic variances of LSEs and MNR estimates of the fundamental frequency in case Model 2 i.i.d. error when sample size $n=100$, $200$, $400$, $500$ and 1000.}   \label{table4}
\begin{tabular}{|c|c|c|c|c|}
\hline
\multicolumn{5}{c}{Sample Size N=100} \\
\hline
$\sigma^2$  & Average & Variance & Asym. Var. (LSE) & Asym. Var. (MNR) \\
\hline
.01   &  .3148      &    9.22e-10   &     1.63e-9   &     4.08e-10  \\ \hline
.25       &   .3148      &    2.34e-8   &     4.08e-8   &     1.02e-8   \\ \hline 
.75       &  .3148      &    7.22e-8   &     1.22e-7   &     3.06e-8    \\ \hline  
1.0       &     .3148       &    9.74e-8   &     1.63e-7   &     4.08e-8  \\ \hline
\multicolumn{5}{c}{Sample Size N=200} \\
\hline
$\sigma^2$  & Average & Variance & Asym. Var. (LSE) & Asym. Var. (MNR) \\
\hline
.01   &     .3143     &    1.31e-10   &     2.04e-10   &     5.10e-11 \\ \hline
.25       &  .3143    &    3.30e-9   &     5.10e-9   &     1.28e-9    \\ \hline 
.75       &      .3143      &    9.75e-9   &     1.53e-8   &     3.82e-9 \\ \hline  
1.0       &    .3143       &    1.28e-8   &     2.04e-8   &     5.10e-9   \\ \hline
\multicolumn{5}{c}{Sample Size N=400} \\
\hline
$\sigma^2$  & Average & Variance & Asym. Var. (LSE) & Asym. Var. (MNR) \\
\hline
.01   &    .3141      &    1.91e-11   &     2.55e-11   &     6.37e-12   \\ \hline
.25       &   .3141     &    4.71e-10   &     6.37e-10   &     1.59e-10   \\ \hline 
.75       &    .3141      &    1.30e-9   &     1.91e-9   &     4.78e-10  \\ \hline  
1.0       &     .3141   &    1.67e-9   &     2.55e-9   &     6.37e-10  \\ \hline
\multicolumn{5}{c}{Sample Size N=500} \\
\hline
$\sigma^2$  & Average & Variance & Asym. Var. (LSE) & Asym. Var. (MNR) \\
\hline
.01   &    .3141      &    1.06e-11   &     1.31e-11   &     3.26e-12 \\ \hline
.25       &   .3141      &    2.55e-10   &     3.26e-10   &     8.16e-11   \\ \hline 
.75       &     .3141       &    6.78e-10   &     9.79e-10   &     2.45e-10 \\ \hline  
1.0       &     .3141       &    8.76e-10   &     1.31e-9   &     3.26e-10   \\ \hline
\multicolumn{5}{c}{Sample Size N=1000} \\
\hline
$\sigma^2$  & Average & Variance & Asym. Var. (LSE) & Asym. Var. (MNR) \\
\hline
.01   &  .3141      &    8.11e-13   &     1.63e-12   &     4.08e-13  \\ \hline
.25       &   .3141      &    2.01e-11   &     4.08e-11   &     1.02e-11   \\ \hline 
.75       &  .3141      &    6.34e-11   &     1.22e-10   &     3.06e-11    \\ \hline  
1.0       &     .3141       &    8.34e-11   &     1.63e-10   &     4.08e-11  \\ \hline
\end{tabular}
\end{table} 
The following are some of the salient features of the numerical experiments reported in Table \ref{table1}-\ref{table4}.
\begin{itemize}
\item[{(i)}] We observe that the average estimators of the fundamental frequency are very close to the true values in all sample sizes and $\sigma^2$ considered.  The estimator has a small positive bias for small sample sizes, but it becomes unbiased when the sample size
is large.
\item[{(ii)}] The variance of the estimate increases with increase in error variance  whereas  it decreases with increase in sample size.  If verifies the consistency property of the proposed estimator.  
\item[{(iii)}] In all the cases considered here, the variance is close to the asymptotic variance of the MNR estimator.  As sample size 
increases it becomes closer.  It is smaller than the asymptotic variance of the LSE in all cases.  Therefore, in line of theorem \ref{theorem1}, improvement is achievable in practice.
\end{itemize}

\section{Data Analysis}  \label{sec5}

In this section, we have analyzed two data sets namely two vowel sounds ``uuu" and ``ahh" using the proposed MNR 
algorithm.  Both the data sets have been obtained from the sound laboratory of the Indian Institute of Technology 
Kanpur.

\subsection{``uuu" data}  \label{sec5.2}

\begin{figure} 
\epsfig{file=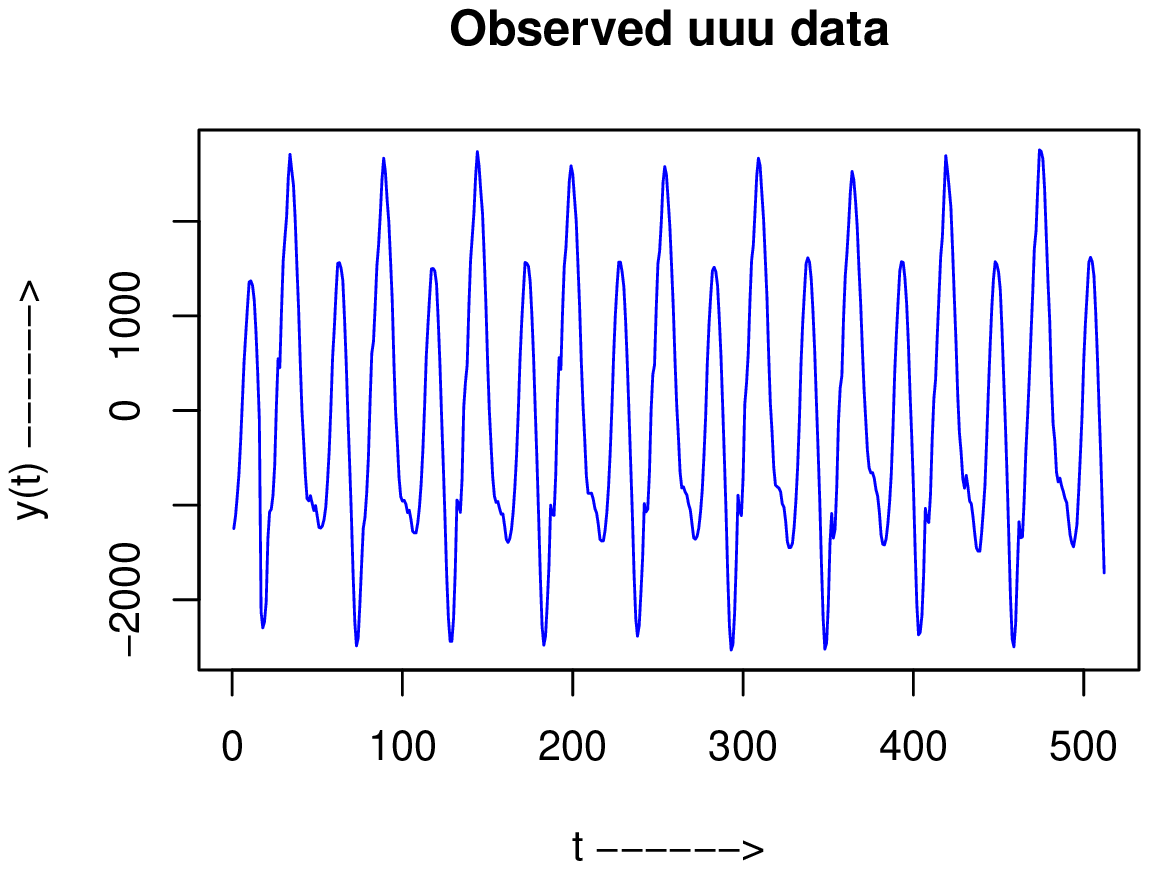, scale=.65} \hspace{.2in}
\epsfig{file=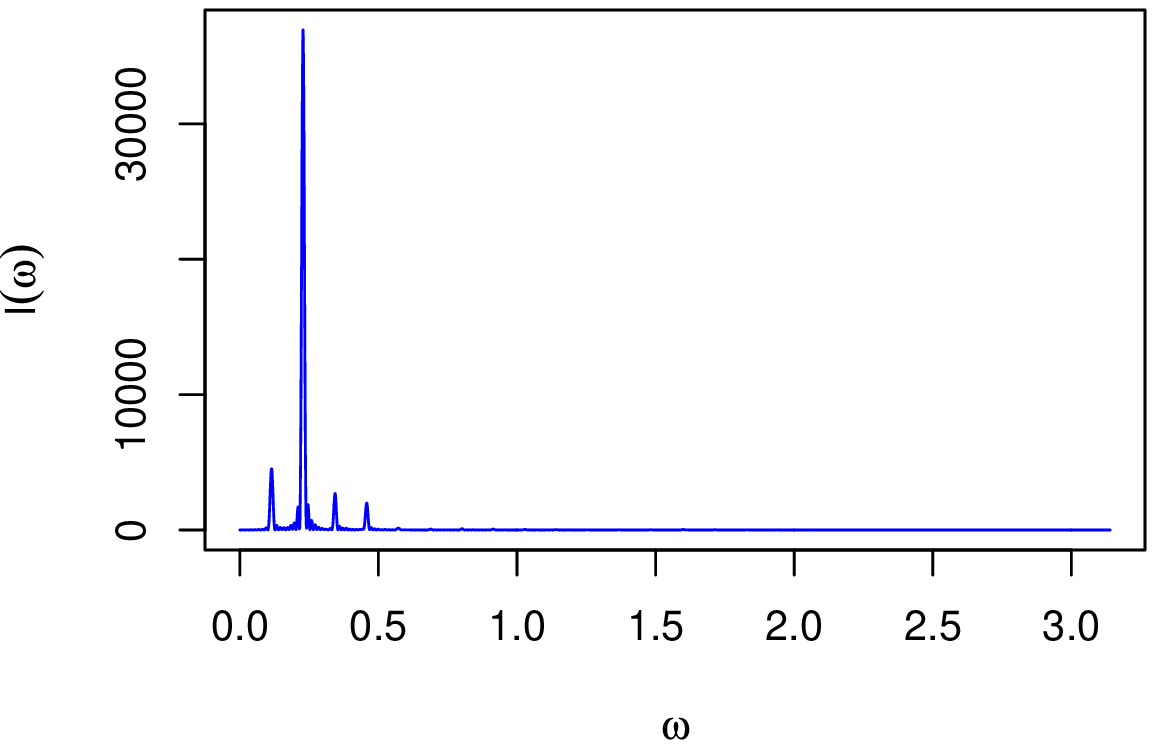, scale=.65}
\caption{Mean corrected ``uuu" data and its periodogram function.}  \label{fig3}
\end{figure}

\begin{figure} 
\epsfig{file=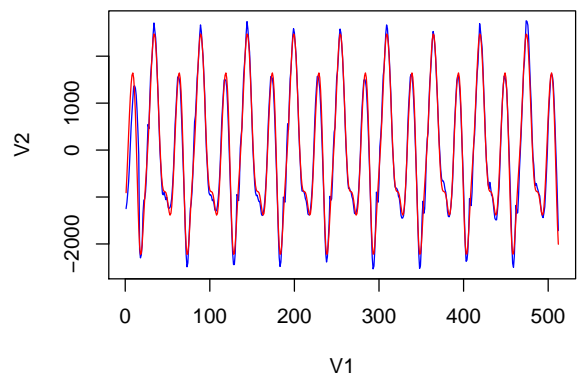, scale=.65}
\caption{The fitted values (red) along with the mean corrected ``uuu" data (blue).}  \label{fig4}
\end{figure}

This data set is for the vowel sound ``uuu".  It contains 512 data points sampled at 10 kHz frequency.  
The mean corrected data and the periodogram function are presented in Fig. \ref{fig3}.  It seems from the periodogram plot that there 
are four harmonics of the fundamental frequency.  Therefore, we take the first significance frequency as the fundamental frequency and 
rest are the harmonics.  Initially, we have estimated the fundamental frequency with four components using the proposed MNR algorithm. 
We check the residuals, but it is not stationary.  After adding two more components sequentially, the error sequence becomes 
stationary.  So, we have fitted the model with $p$ = 6.  The MNR estimate of $\lambda$ is $0.1142$.  The linear parameters are estimated as mentioned before.  The fitted values (red)and mean corrected observed ``uuu" data (blue) are plotted in Fig. \ref{fig4}.  They match very well.    Using the parameter estimates, the error sequence is estimated which can be fitted  as the following stationary ARMA(2,4) process;
\beanno
e(t) &=& -4.636 + 1.8793 e(t-1) -0.9308 e(t-2) + \epsilon(t) -0.8657 \epsilon(t-1) \\
      &  &- 0.4945 \epsilon(t-2) + 0.2859 \epsilon(t-3) +0.1415 \epsilon(t-4).
\eeanno
\subsection{``ahh" data} \label{sec5.3}
\begin{figure} 
\epsfig{file=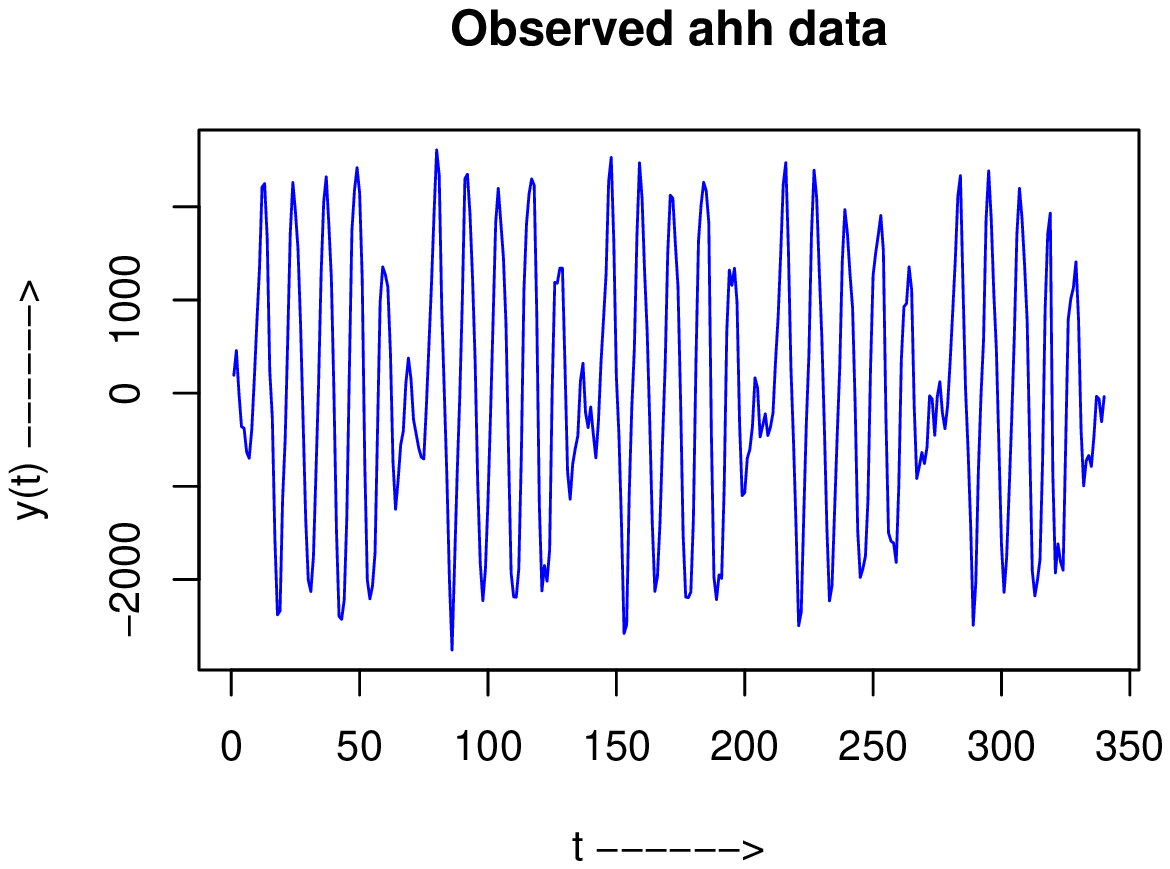, scale=.65} \hspace{.2in}
\epsfig{file=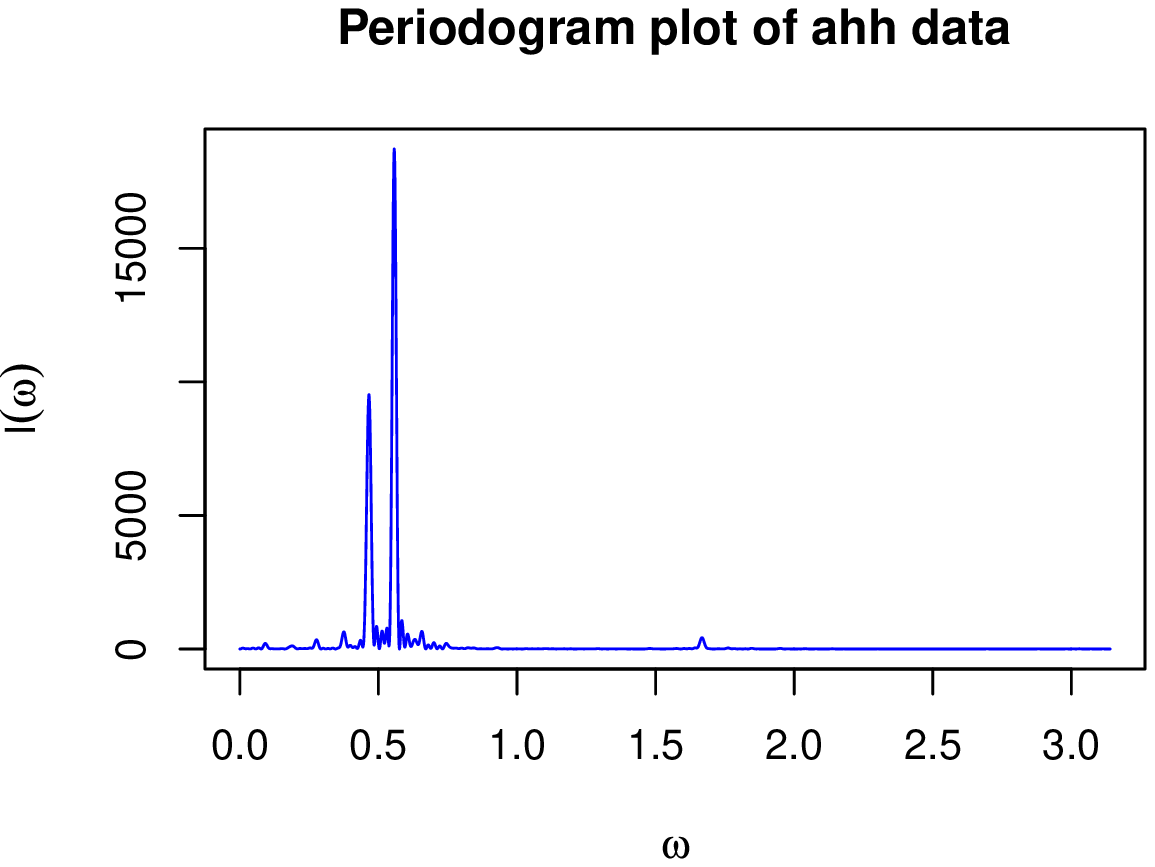, scale=.65}
\caption{Mean corrected ``ahh" data and its periodogram function.}  \label{fig5}
\end{figure}

This is a sound data ``ahh".  It contains 340 signal values sampled at 10 kHz frequency,  The mean corrected data and its periodogram function are plotted in Fig. \ref{fig5}.  Following the same methodology as applied in case of ``uuu" data, it is observed that the fundamental 
frequency model with $p$ = 6, fits the data quite well and the error sequence also becomes stationary.
Based on the MNR algorithm, the estimate of $\lambda$ is obtained as $.0929$.  Then, the fitted values are obtained similarly as ``uuu" 
data set.  The fitted values match quite well with the mean corrected ``ahh" data.  The estimated error in this case is 
$$
e(t) = 1.8128 + 0.6816 e(t-1) + \epsilon(t) +0.4246 \epsilon(t-1) -0.5315 \epsilon(t-2) - 0.6572 \epsilon(t-3).
$$ 
This is a stationary ARMA(1,3) process and  can be expressed as \eqref{assum_eq}.
\begin{figure} 
\epsfig{file=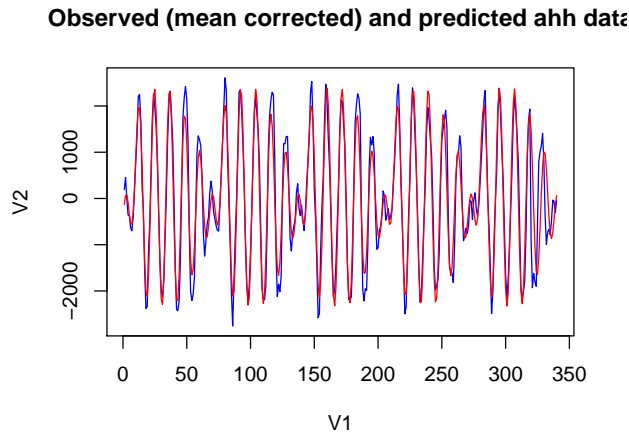, scale=.65}
\caption{The fitted values (red) along with the mean corrected ``ahh" data (blue).}  \label{fig6}
\end{figure}

\section{Concluding Remarks}  \label{sec6}
In this paper, we have considered the fundamental frequency model.  This model is the multiple sinusoidal frequency model, where frequencies are harmonics of a fundamental frequency.  We are mainly interested in estimating $\lambda$, the fundamental frequency.  Once $\lambda$ is estimated efficiently, the other linear parameters are easily estimated using LS or approximate LS approach.  It is well known that the 
NR algorithm does not work well in case of the sinusoidal model.  In this paper, we propose to modify the step factor in NR algorithm and observe that it improves the performance of the algorithm quite effectively.  The asymptotic variance of the proposed estimator is smaller than 
the asymptotic variance of the LSE.    The fundamental frequency as a single nonlinear parameter has a quite complicated form in LS or approximate LS approach.  The modified NR algorithm does not require any optimization. The calculation of first and second order derivatives at each step is only required, hence it is very simple to implement.  

We think sequential application of the proposed algorithm will be required  if higher order harmonic terms are present for more than one fundamental frequency.  The proposed algorithm can be extended in case of multiple fundamental frequency model (Chirstensen et  al. \cite{CHJJ:2011}) and cluster type model (Nandi and Kundu \cite{NK:2006a}).   This a topic of ongoing research and would be reported elsewhere. 

\noindent Acknowledgements: The authors would like to thank the reviewers and the Associate Editor for their constructive 
comments which have helped to improve the manuscript significantly.

\section*{Appendix}

\noindent {\bf Proof of Theorem \ref{theorem1}}:  In the proof of theorem \ref{theorem1}, at any iteration we use $\wt{\lambda}$ as the initial estimator and $\wh{\lambda}$ as the updated estimator of $\lambda$.   Now, define the following matrices to express the first and second order derivatives of $R_j(\lambda)$. 
$$
{\bf D}_j = \hbox{diag}\{j, 2j, \ldots, nj\}, \ \ \ {\bf E} = \left [ \begin{matrix}0 & 1 \cr -1 & 0 \cr \end{matrix}
\right ], \ \ \ \stackrel{.}{\bf X}_j = \frac{d}{d\lambda} {\bf X} = {\bf D}_j {\bf X}_j {\bf E}, \ \ \  
\stackrel{..}{\bf X}_j = \frac{d^2}{d\lambda^2} {\bf X}= -{\bf D}_j^2 {\bf X}_j.
$$
Note that $\ds {\bf EE }= -{\bf I}$, $\ds {\bf EE}^T = {\bf I} = {\bf E}^T{\bf E}$ and 
$$
\frac{d}{d\lambda} ({\bf X}_j^T{\bf X}_j)^{-1} = -({\bf X}_j^T{\bf X}_j)^{-1}[{\bf \dotx}_j^T{\bf X}_j + {\bf X}_j^T{\bf \dotx}_j]({\bf X}_j^T{\bf X}_j)^{-1}.
$$
Write $\ds \frac{d}{d \lambda} R_j(\lambda) = R'_j(\lambda)$ and $\frac{d^2}{d^2 \lambda} R_j(\lambda) = R''_j(\lambda)$.  Then 
\begin{equation}
\frac{1}{2} R'_j(\lambda) =  {\bf Y}^T {\bf \dotx}_j ({\bf X}_j^T{\bf X}_j)^{-1} {\bf X}^T{\bf Y} - {\bf Y}^T{\bf X}_j({\bf X}_j^T{\bf X}_j)^{-1}{\bf \dotx}_j^T{\bf X}_j({\bf X}_j^T{\bf X}_j)^{-1}{\bf X}_j^T{\bf Y},   \label{rprime_eq}
\end{equation}
and
\bea
\frac{1}{2} \ R''_j(\lambda) & = &  {\bf Y}^T {\bf \ddotx}_j ({\bf X}_j^T{\bf X}_j)^{-1} {\bf X}_j^T{\bf Y} - 
 {\bf Y}^T{\bf \dotx}_j({\bf X}_j^T{\bf X}_j)^{-1}({\bf \dotx}_j^T{\bf X}_j + {\bf X}_j^T{\bf \dotx}_j)({\bf 
X_j}^T{\bf X}_j)^{-1}{\bf X}_j^T{\bf Y}    \nonumber  \\
& + &  {\bf Y}^T{\bf \dotx}_j({\bf X}_j^T{\bf X}_j)^{-1}{\bf \dotx}_j^T{\bf Y} - {\bf Y}^T{\bf \dotx}_j({\bf X}_j^T{\bf X}_j)^{-1}{\bf \dotx}_j^T{\bf X}_j({\bf X}_j^T{\bf X}_j)^{-1} {\bf X}_j^T{\bf Y}  \nonumber  \\
& + & {\bf Y}^T{\bf X}_j({\bf X}_j^T{\bf X}_j)^{-1}({\bf \dotx}_j^T{\bf X}_j + {\bf X}_j^T{\bf \dotx}_j)({\bf X}_j^T{\bf X}_j)^{-1}{\bf \dotx}_j^T{\bf X}_j({\bf X}_j^T{\bf X}_j)^{-1} {\bf X}_j^T{\bf Y}   \nonumber  \\
& - &  {\bf Y}^T{\bf X}_j({\bf X}_j^T{\bf X}_j)^{-1}({\bf \ddotx_j}^T{\bf X}_j)({\bf X}_j^T{\bf X}_j)^{-1} {\bf X}_j^T{\bf Y}    \nonumber  \\ 
&- & {\bf Y}^T{\bf X}_j({\bf X}_j^T{\bf X}_j)^{-1}({\bf \dotx}_j^T{\bf \dotx}_j)({\bf X}_j^T{\bf X}_j)^{-1} {\bf X}_j^T{\bf Y}    \nonumber  \\
& + &{\bf Y}^T{\bf X}_j({\bf X}_j^T{\bf X}_j)^{-1}{\bf \dotx}_j^T{\bf X}_j({\bf X}_j^T{\bf X}_j)^{-1}({\bf \dotx}_j^T{\bf X}_j + {\bf X}_j^T{\bf \dotx}_j)({\bf X}_j^T{\bf X}_j)^{-1} {\bf X}_j^T{\bf Y}   \nonumber  \\
& - & {\bf Y}^T{\bf X}_j({\bf X}_j^T{\bf X}_j)^{-1}{\bf \dotx}_j^T{\bf X}_j({\bf X}_j^T{\bf X}_j)^{-1} {\bf \dotx}_j^T{\bf Y}.  \label{rprimeprime_eq}
\eea
By definition $\ds g(\lambda) =  \sum_{j=1}^p R_j(\lambda)$, therefore,  we have $\ds g'(\lambda) = \sum_{j=1}^p R'_j(\lambda)$ and $\ds g''(\lambda) = \sum_{j=1}^p R''_j(\lambda)$. 

\noindent Assume that $\ds \widetilde{\lambda} - \lambda = O_p (n^{-1 - \delta})$, $\ds \delta \in (0, \frac{1}{2}]$.   Therefore, for large $n$,  at $\lambda = \widetilde{\lambda}$, 
\begin{equation}
\displaystyle (\frac{1}{n} {\bf X}_j^T{\bf X}_j)^{-1} = (\frac{1}{n} {\bf X}_j(\widetilde{\lambda})^T{\bf X}_j(\widetilde{\lambda}))^{-1} = 2 \ I + O_p(\frac{1}{n}).  \label{proof_eq1}
\end{equation}
Using the large sample approximation \eqref{proof_eq1} in the first term of $\frac{1}{2} R'_j(\lambda)$ in \eqref{rprime_eq}, we have at $\lambda = \widetilde{\lambda}$,   
\beanno
&&  \frac{1}{n^3} {\bf Y}^T {\bf \dotx}_j ({\bf X}_j^T{\bf X}_j)^{-1} {\bf X}_j^T{\bf Y} \\ &=&   
\frac{1}{n^3} {\bf Y}^T {\bf \dotx}_j (\wtilde{\lambda}) ({\bf X}_j (\wtilde{\lambda})^T{\bf X}_j(\wtilde{\lambda}))^{-1} {\bf X}_j(\wtilde{\lambda})^T{\bf Y} \\
&=& \frac{2}{n^4}{\bf Y}^T{\bf D}_j{\bf X}_j(\wtilde{\lambda}) {\bf E X}_j(\wtilde{\lambda})^T{\bf Y} \\
& = & \frac{2j}{n^4} \left[\left ( \sum_{t=1}^n y(t) t \cos(\widetilde{j \lambda} t) \right )
\left ( \sum_{t=1}^n y(t) \sin(j \widetilde{\lambda} t) \right ) 
-  \left ( \sum_{t=1}^n y(t) t \sin(j \widetilde{\lambda} t) \right )
\left ( \sum_{t=1}^n y(t) \cos(j \widetilde{\lambda} t) \right ) \right].
\eeanno
Then along the same line as Kundu {\it et al.} \cite{KBNB:2011}, it can be shown that, 
\be
\sum_{t=1}^n y(t) \cos(j \widetilde{\lambda} t) = \frac{n}{2} \left ( A_j + O_p(n^{-\delta}) \right ), \ \ \ \ 
\sum_{t=1}^n y(t) \sin(j \widetilde{\lambda} t) = \frac{n}{2} \left ( B_j + O_p(n^{-\delta}) \right ).  \label{proof_eq11}
\ee
Now consider
\bea
\sum_{t=1}^n y(t) t e^{-i j \widetilde{\lambda} t} & = & \sum_{t=1}^n \left ( \sum_{k=1}^p \left[A_k 
\cos(k \lambda t) + B_k \sin (k \lambda t) + e(t)\right] \right ) t e^{-i j \widetilde{\lambda} t}     \nonumber \\
& = & \frac{1}{2} \sum_{k=1}^p (A_k - iB_k) \sum_{t=1}^n t \ e^{i (k\lambda - j \widetilde{\lambda}) t} +   \nonumber \\
&& \frac{1}{2} \sum_{k=1}^p  (A_k + iB_k) \sum_{t=1}^n t \ e^{-i (k\lambda + j \widetilde{\lambda}) t} + 
\sum_{t=1}^n e(t) t e^{-i j \widetilde{\lambda} t}  \label{proof_eq2}
\eea
Similarly as Bai {\it et al.}\cite{BRCK:2003}, the following can be established for harmonics of fundamental frequency;
\beanno
&&\sum_{t=1}^n t \ e^{-i (k \lambda + j \widetilde{\lambda}) t} = O_p(n),   ~~~\forall ~~k,j=1,\ldots p \\
&&\sum_{t=1}^n t \ e^{-i (k \lambda - j \widetilde{\lambda}) t} = O_p(n),   ~~~\forall ~~k \ne j=1,\ldots p
\eeanno
and for $\ds k=j$,
\bea
\sum_{t=1}^n t \ e^{i (\lambda - \widetilde{\lambda})jt} &=& \sum_{t=1}^n t + i (\lambda - \widetilde{\lambda}) j \sum_{t=1}^n t^2 - \frac{1}{2}(\lambda - \widetilde{\lambda})^2 j^2 \sum_{t=1}^n t^3      \nonumber  \\
& &  - \frac{1}{6} i (\lambda - \widetilde{\lambda})^3 j^3 \sum_{t=1}^n t^4  +  
\frac{1}{24} (\lambda - \widetilde{\lambda})^4 j^4 \sum_{t=1}^n t^5 e^{i (\lambda - \lambda^*)j t}.  \label{proof_eq3}
\eea
The last term of \eqref{proof_eq3} is approximated as $$ \frac{1}{24} (\lambda - \widetilde{\lambda})^4 j^4 \sum_{t=1}^n t^5 e^{i (\lambda - \lambda^*)j t} =  O_p(n^{2-4\delta}).$$  

For the last term in \eqref{proof_eq2},  choose $L$ large enough such that $L \delta > 1$ and using the Taylor series expansion of
$e^{-i j \widetilde{\lambda} t}$ we obtain,
\beanno
& &\sum_{t=1}^n e(t) t e^{-i j \widetilde{\lambda} t} \\
& =  & \sum_{m=0}^{\infty} a(m) \sum_{t=1}^n e(t-m) t 
e^{-i j \widetilde{\lambda} t}  \\
& =  & \sum_{m=0}^{\infty} a(m) \sum_{t=1}^n e(t-m) t e^{-i j \lambda t} + 
\sum_{m=0}^{\infty} a(m) \sum_{l=1}^{L-1} \frac{(-i(\widetilde{\lambda} - \lambda)j)^l}{l!}
\sum_{t=1}^n e(t-m) t^{l+1} e^{-i j \lambda t}  \\
&  & + \sum_{m=0}^{\infty} a(m) \frac{\theta (n(\widetilde{\lambda} - \lambda))^L}{L!} 
\sum_{t=1}^n t |e(t-m)| \  (\hbox{here $|\theta| < 1$})  \\
& = & \sum_{m=0}^{\infty} a(m) \sum_{t=1}^n e(t-m) t e^{-i j\lambda t} + 
\sum_{l=1}^{L-1} O_p(n^{-(1+\delta) l}) O_p(n^{l+\frac{3}{2}}) + 
\sum_{m=0}^{\infty} a(m) O_p(n^{\frac{5}{2} - L \delta})  \\
& = & \sum_{m=0}^{\infty} a(m) \sum_{t=1}^n e(t-m) t e^{-i j \lambda t} +  
O_p(n^{\frac{5}{2} - L \delta}).
\eeanno
Therefore,
\bea
&&\sum_{t=1}^n y(t) t \cos(\widetilde{j \lambda} t) \nonumber \\
& = & \frac{1}{2} \left [\sum_{k=1}^p  
A_k \left ( \sum_{t=1}^n t - \frac{1}{2}(\lambda - \widetilde{\lambda})^2 j^2 \sum_{t=1}^n t^3 \right ) \right.  \nonumber   \\ 
& & \left.+ \sum_{k=1}^p B_k \left ( \sum_{t=1}^n (\lambda - \widetilde{\lambda}) j t^2 - \frac{1}{6} 
(\lambda - \widetilde{\lambda})^3 j^3 \sum_{t=1}^n t^4 \right ) \right ]   \nonumber  \\
& & + \sum_{m=0}^{\infty} a(m) \sum_{t=1}^n e(t-m) t \cos(j \lambda t) 
+  O_p(n^{\frac{5}{2} - L \delta}) + O_p(n) + O_p(n^{2-4 \delta}).   \label{proof_eq31}
 \eea
Similarly,
\bea
&&\sum_{t=1}^n y(t) t \sin(\widetilde{j \lambda} t)  \nonumber \\
 & = & \frac{1}{2} \left [ \sum_{k=1}^p 
B_k \left ( \sum_{t=1}^n t - \frac{1}{2}(\lambda - \widetilde{\lambda})^2 j^2  \sum_{t=1}^n t^3 \right ) \right.   \nonumber   \\
&& \left.-\sum_{k=1}^p A_k \left ( \sum_{t=1}^n (\lambda - \widetilde{\lambda}) j t^2 - \frac{1}{6} 
(\lambda - \widetilde{\lambda})^3 j^3 \sum_{t=1}^n t^4 \right ) \right ]    \nonumber   \\
& & + \sum_{m=0}^{\infty} a(m) \sum_{t=1}^n e(t-m) t \sin(j \lambda t) 
+  O_p(n^{\frac{5}{2} - L \delta}) + O_p(n) + O_p(n^{2-4 \delta}).   \label{proof_eq32}
 \eea
Next, the second term of $\frac{1}{2} R'_j(\lambda)$ in \eqref{rprime_eq} is approximated as
\bea
& & \frac{1}{n^3} {\bf Y}^T{\bf X}_j({\bf X}_j^T{\bf X}_j)^{-1}{\bf \dotx}_j^T{\bf X}_j({\bf X}_j^T{\bf X}_j)^{-1}{\bf X}_j^T{\bf Y} \nonumber  \\ 
& = &  \frac{1}{n^3} {\bf Y}^T{\bf X}_j({\bf X}_j^T{\bf X}_j)^{-1}{\bf E}^T{\bf X}_j^T{\bf D}_j{\bf X}_j({\bf X}_j^T{\bf X}_j)^{-1}{\bf X}_j^T{\bf Y}  \nonumber \\
& = & \frac{1}{n^3}{\bf Y}^T{\bf X}_j \left(2 {\bf I} + O_p(\frac{1}{n}) \right){\bf E}^T \left(\frac{1}{4}{\bf I} + O_p(\frac{1}{n}) \right) \left(2 {\bf I} + O_p(\frac{1}{n})\right){\bf X}_j^T{\bf Y}  \nonumber \\
& = & \frac{j}{n^3}{\bf Y}^T{\bf X}_j{\bf E}^T{\bf X}_j^T{\bf Y} + O_p(\frac{1}{n}) = O_p(\frac{1}{n}),  \label{proof_eq4}
\eea
for large $n$ and $\lambda=\wtilde{\lambda}$.


Now to simplify  $R'_j(\widetilde{\lambda})$ and $R''_j(\widetilde{\lambda})$, we need the following results, for any $\lambda \in (0,\pi)$.
\bea
 &&\sum_{t=1}^n t \ \cos^2(j\lambda t) = \frac{n^2}{4} + O(n), \ \ \ \ 
\sum_{t=1}^n t \ \sin^2(j\lambda t) = \frac{n^2}{4} + O(n),  \label{results1}
\\
 &&\sum_{t=1}^n \cos^2(j\lambda t) = \frac{n}{2} + o(n), \ \  \ \ 
\sum_{t=1}^n \sin^2(j\lambda t) = \frac{n}{2} + o(n), \label{results2}
\\
 &&\sum_{t=1}^n t^2 \ \cos^2(j \lambda t) = \frac{n^3}{6} + O(n^2), \ \  \ \
\sum_{t=1}^n t^2 \ \sin^2(j \lambda t) = \frac{n^3}{6} + O(n^2), \label{results3}
\eea
and
\bea
 &&\frac{1}{n^2} {\bf Y}^T {\bf D}_j {\bf X}_j = \frac{j}{4}(A_j \ \ B_j) + O_p(\frac{1}{n}), \ \ \ 
 \frac{1}{n^3}  {\bf Y}^T {\bf D}_j^2 {\bf X}_j = \frac{j^2}{6}(A_j \ \ B_j) + O_p(\frac{1}{n}),  \label{results4}
\\
&& \frac{1}{n^3}{\bf X}_j^T{\bf D}_j^2{\bf X}_j = \frac{j^2}{6} {\bf I} + O_p(\frac{1}{n}), \ \ \ \ \frac{1}{n}{\bf X}_j^T{\bf Y} = \frac{1}{2}(A_j \ \ B_j)^T + O_p(\frac{1}{n}),  \label{results5}
\\
&& \frac{1}{n^2}{\bf X}_j^T{\bf D}_j {\bf X}_j = \frac{j}{4} {\bf I} + 
O_p(\frac{1}{n}). \label{results6}
\eea

Next to simplify $\frac{1}{2n^3} \ R''_j(\widetilde{\lambda}) $, use \eqref{proof_eq1} at the first step.
\beanno
\frac{1}{2n^3} \ R''_j(\widetilde{\lambda}) & = & \frac{2}{n^4} {\bf Y}^T{\bf \ddotx}_j {\bf X}_j^T{\bf Y} - 
\frac{4}{n^5} {\bf Y}^T{\bf \dotx}_j({\bf \dotx}_j^T{\bf X}_j + {\bf X}_j^T{\bf \dotx}_j){\bf X}_j^T{\bf Y}_j + \frac{2}{n^4} {\bf Y}^T{\bf \dotx}_j {\bf \dotx}_j^T{\bf Y}_j  \\
& - & \frac{4}{n^5} {\bf Y}^T{\bf \dotx}_j {\bf \dotx}_j^T{\bf X}_j {\bf X}_j^T{\bf Y} + \frac{8}{n^6} 
{\bf Y}^T{\bf X}_j({\bf \dotx}_j^T{\bf X}_j + {\bf X}_j^T{\bf \dotx}_j){\bf \dotx}_j^T{\bf X}_j {\bf X}_j^T{\bf Y}   \\
& - & \frac{4}{n^5} {\bf Y}^T{\bf X}_j{\bf \ddotx}_j^T{\bf X}_j {\bf X}_j^T{\bf Y} - \frac{4}{n^5} {\bf Y}^T{\bf X}_j {\bf \dotx}_j^T{\bf \dotx}_j 
{\bf X}_j^T{\bf Y} \\ 
&+ &\frac{8}{n^6}{\bf Y}^T{\bf X}_j {\bf \dotx}_j^T{\bf X}_j({\bf \dotx}_j^T{\bf X}_j + {\bf X}_j^T{\bf \dotx}_j){\bf X}_j^T{\bf Y}  
-  \frac{4}{n^5}  {\bf Y}^T{\bf X}_j {\bf \dotx}_j^T{\bf X}_j {\bf \dotx}_j^T{\bf Y} + O_p(\frac{1}{n}).  
\eeanno
In the second step, use ${\bf \ddotx}_j = {\bf D}_j {\bf X}_j {\bf E}$ and ${\bf \ddotx}=-{\bf D}_j^2 {\bf X}_j$.
\beanno
\frac{1}{2n^3} \ R''_j(\widetilde{\lambda}) & = &  - \frac{2}{n^4}  {\bf Y}^T {\bf D}_j^2 {\bf X}_j {\bf X}_j^T{\bf Y} - 
\frac{4}{n^5} {\bf Y}^T {\bf D}_j {\bf X}_j {\bf E} ({\bf E}^T{\bf X}_j^T{\bf D}_j {\bf X}_j + {\bf X}_j^T{\bf D}_j {\bf X}_j {\bf E}){\bf X}_j^T{\bf Y} \\
& + & \frac{2}{n^4} {\bf Y}^T{\bf D}_j {\bf X}_j {\bf EE}^T{\bf X}_j^T{\bf D}_j {\bf Y}  -  \frac{4}{n^5} {\bf Y}^T{\bf D}_j {\bf X}_j {\bf EE}^T{\bf X}_j^T{\bf D}_j {\bf X}_j {\bf X}_j^T{\bf Y}  \\
& + & \frac{8}{n^6} {\bf Y}^T{\bf X}_j({\bf E}^T{\bf X}_j^T{\bf D}_j {\bf X}_j + {\bf X}_j^T{\bf D}_j {\bf X}_j {\bf E}){\bf E}^T{\bf X}_j^T
{\bf D}_j {\bf X}_j {\bf X}_j^T{\bf Y}  + \frac{4}{n^5} {\bf Y}^T{\bf X}_j {\bf X}_j^T{\bf D}_j^2{\bf X}_j {\bf X}_j^T{\bf Y}  \\
& - & \frac{4}{n^5} {\bf Y}^T{\bf X}_j {\bf E}^T{\bf X}_j^T{\bf D}_j^2{\bf X}_j {\bf E X}_j^T{\bf Y} + \frac{8}{n^6} {\bf Y}^T{\bf X}_j {\bf E}^T{\bf X}_j^T{\bf D}_j {\bf X}_j({\bf E}^T{\bf X}_j^T{\bf D}_j {\bf X}_j \\ 
&+ & {\bf X}_j^T{\bf D}_j {\bf X}_j {\bf E}){\bf X}_j^T{\bf Y} 
 -  \frac{4}{n^5}  {\bf Y}^T{\bf X}_j {\bf E}^T{\bf X}_j^T{\bf D}_j {\bf X}_j {\bf E}^T{\bf X}_j^T{\bf D}_j {\bf Y} +  O_p(\frac{1}{n}).
\eeanno
Next, using \eqref{results1}-\eqref{results6}, we observe
\bea
\frac{1}{2n^3} \ R''_j(\widetilde{\lambda}) &=& (A_j^2+B_j^2) \left [ - \frac{j^2}{6} - 0 + \frac{j^2}{8} - 
\frac{j^2}{8} + 0 + \frac{j^2}{6} - \frac{j^2}{6} + 0 + \frac{j^2}{8}
\right ] + O_p(\frac{1}{n})   \nonumber \\
&=& -\frac{j^2}{24}(A_j^2+B_j^2) + O_p(\frac{1}{n}).  \label{proof_eq6}
\eea
The correction factor in Newton-Raphson algorithm can be written as 
\be
 \frac{g'(\wt{\lambda})}{g''(\wt{\lambda})} = \frac{\ds \frac{1}{2n^3} \sum_{j=1}^p R'_j(\wt{\lambda})}{\ds \frac{1}{2n^3} \sum_{j=1}^p R''_j(\wt{\lambda})}       \label{gratio_eq}
\ee
Using \eqref{proof_eq31}, \eqref{proof_eq32} and \eqref{proof_eq4}, $\ds \frac{1}{2n^3} \sum_{j1}^p R'_j(\wt{\lambda})$ is simplified as 
\beanno
\frac{1}{2n^3} \sum_{j=1}^p R'_j(\wt{\lambda}) &=&
 \frac{2}{n^4} \sum_{j=1}^p j\left[\frac{n}{2}(B_j +O_p(n^{-\delta})) \left\{\frac{A_j}{2}\left(\sum_{t=1}^n t - \frac{1}{2}(\lambda - \widetilde{\lambda})^2 j^2 \sum_{t=1}^n t^3 \right) \right. \right. \\
&&    + \frac{B_j}{2} \left(\sum_{t=1}^n (\lambda - \widetilde{\lambda}) j t^2 - \frac{1}{6} (\lambda - \widetilde{\lambda})^3 j^3 \sum_{t=1}^n t^4 \right) \\ 
&& + \sum_{k=0}^\infty a(k) \sum_{t=1}^n e(t-k) t \cos(j \lambda t) + O_p(n^{\frac{5}{2} - L\delta})  + O_p(n) + O_p(n^{2-4\delta}) \BRB \\
&& - \frac{n}{2}(A_j +O_p(n^{-\delta})) \left\{ \frac{B_j}{2}\left(\sum_{t=1}^n t - \frac{1}{2}(\lambda - \widetilde{\lambda})^2 j^2 \sum_{t=1}^n t^3 \right) \right. \\
& &  - \frac{A_j}{2} \left(\sum_{t=1}^n (\lambda - \widetilde{\lambda}) j t^2 - \frac{1}{6} (\lambda - \widetilde{\lambda})^3 j^3 \sum_{t=1}^n t^4 \right) 
\\
&& + \sum_{k=-}^\infty a(k) \sum_{t=1}^n e(t-k) t \sin(j \lambda t) + O_p(n^{\frac{5}{2} - L\delta})   + O_p(n) + O_p(n^{2-4\delta}) \BRB \BRS \\
&=& \sum_{j=1}^p j \left[\frac{1}{2 n^3} (A_j^2 + B_j^2) \left\{ \sum_{t=1}^n (\lambda - \widetilde{\lambda}) j t^2 -\frac{1}{6} (\lambda - \widetilde{\lambda})^3 j^3 \sum_{t=1}^n t^4 \right\}  \right. 
\\
& & + \frac{1}{ n^3} \left\{ B_j \sum_{k=0}^\infty a(k) \sum_{t=1}^n e(t-k) t \cos(j \lambda t) \right. \\
&& \left. \left.+ A_j \sum_{k=0}^\infty a(k) \sum_{t=1}^n e(t-k) t \sin(j \lambda t) \right\} \right] \\ 
& & + O_p(n^{-\frac{1}{2} - L \delta}) + O_p(n^{-2}) + O_p(n^{-1 -4\delta}),
\eeanno
and using \eqref{proof_eq6}, the denominator of \eqref{gratio_eq} is 
$\ds 
\frac{1}{2n^3} \sum_{j=1}^p R''_j(\wt{\lambda}) = - \frac{1}{24} \sum_{j=1}^p j^2 (A_j^2 + B_j^2)  +  O_p(\frac{1}{n})$. 
Therefore,
\bea
\wh{\lambda} &=& \wt{\lambda} - \frac{1}{4}\frac{g'(\wt{\lambda})}{g''(\wt{\lambda})} 
= \wt{\lambda} - \frac{1}{4} \frac{\ds \frac{1}{2n^3} \sum_{j=1}^p R'_j(\wt{\lambda})}{\ds \frac{1}{2n^3} \sum_{j=1}^p R''_j(\wt{\lambda})}   \hspace{3in} \nonumber \\
&=& \wt{\lambda} - \frac{1}{4} \frac{\ds \frac{1}{2n^3} \sum_{j=1}^p R'_j(\wt{\lambda})}{\ds  -\frac{1}{24}\sum_{j=1}^p j^2 (A_j^2 + B_j^2) +O_p(\frac{1}{n})}    \nonumber \\
&=& \wt{\lambda} + 6 \frac{\ds \frac{1}{2n^3} \sum_{j=1}^p R'_j(\wt{\lambda})}{\ds  \beta^* +O_p(\frac{1}{n})}    \nonumber \\
&=& \wt{\lambda} + \frac{6}{(\beta^*+ O_p(\frac{1}{n}))} \sum_{j=1}^p j \left[\frac{1}{2 n^3} (A_j^2 + B_j^2) \left\{ \sum_{t=1}^n (\lambda - \widetilde{\lambda}) j t^2 -\frac{1}{6} (\lambda - \widetilde{\lambda})^3 j^3 \sum_{t=1}^n t^4 \right\}  \right]  \nonumber \\
&& \hspace{.3in }+ \frac{6}{(\beta^* + O_p(\frac{1}{n}))} \sum_{j=1}^p j \frac{1}{ n^3} \left\{ B_j \sum_{k=0}^\infty a(k) \sum_{t=1}^n e(t-k) t \cos(j \lambda t) \right.   \nonumber  \\
&&  \hspace{.3in} \left.+ A_j \sum_{k=0}^\infty a(k) \sum_{t=1}^n e(t-k) t \sin(j \lambda t) \right\}+ O_p(n^{-\frac{1}{2} - L \delta}) + O_p(n^{-2}) + O_p(n^{-1 -4\delta}).   \nonumber  \\
&=& \lambda + (\lambda - \wt{\lambda}) O_p(n^{-2 \delta})  \nonumber \\
&& \hspace{.3in }+ \frac{6}{(\beta^* + O_p(\frac{1}{n}))} \sum_{j=1}^p j \frac{1}{ n^3} \left\{ B_j \sum_{k=0}^\infty a(k) \sum_{t=1}^n e(t-k) t \cos(j \lambda t) \right.   \nonumber  \\
&&  \hspace{.3in} \left.+ A_j \sum_{k=0}^\infty a(k) \sum_{t=1}^n e(t-k) t \sin(j \lambda t) \right\}+ O_p(n^{-\frac{1}{2} - L \delta}) + O_p(n^{-2}) + O_p(n^{-1 -4\delta}). \nonumber \\
  \label{proof_eq7}
\eea
Here $\beta^* = \ds \sum_{j=1}^p j^2 (A_j^2 + B_j^2)$ is same as defined after \eqref{asym_lse}.  When $\delta \le \frac{1}{6}$ in \eqref{proof_eq7},  $\wh{\lambda} - \lambda = O_p(n^{-1-3\delta})$ whereas if $\delta > \frac{1}{6}$, then for large $n$,
\beanno
n^{3/2} (\wh{\lambda} - \lambda) &\stackrel{d}{=}&  \frac{6 n^{-3/2}}{\beta^*} \sum_{j=1}^p j\left\{ B_j \sum_{k=0}^\infty a(k) \sum_{t=1}^n e(t-k) t \cos(j \lambda t) \right. \\
&& \hspace{1in} \left. + A_j \sum_{k=0}^\infty a(k) \sum_{t=1}^n e(t-k) t \sin(j \lambda t) \right\} \\
\stackrel{d}{\longrightarrow} \mathcal{N}(0, \gamma)
\eeanno
where 
\beanno
\gamma &=& \frac{36}{{\beta^*}^2} \frac{\sigma^2}{6} \sum_{j=1}^p  j^2 (A_j^2 + B_j^2) \left[\left\{ \sum_{k=0}^\infty a(k) \cos(kj\lambda) \right\}^2 + \left\{\sum_{k=0}^\infty a(k) \sin(kj\lambda) \right\}^2 \right] \\
&=& \frac{6}{{\beta^*}^2} \sigma^2  \sum_{j=1}^p  j^2 (A_j^2 + B_j^2) \left| \sum_{k=0}^\infty a(k) e^{-ikj\lambda} \right|^2 = \frac{6\sigma^2 \delta_G}{{\beta_*}^2}.
\eeanno
This proves the theorem.

\end{document}